\documentclass[a4paper]{article}
\usepackage{amssymb}
\usepackage{amsmath}

\input tcilatex

\begin{document}

\author{Fang Li \\
Department of Mathematics, Jilin University \and Haibao Duan\thanks{%
Supported by NSFC} \\
Institute of Mathematics, Chinese Academy of Sciences}
\title{Homology rigidity of Grassmannians\\
{\small Dedicated to Professor Wenjun Wu on his 90th birthday}}
\date{ }
\maketitle

\begin{abstract}
Applying the theory of Gr\"{o}bner basis to the Schubert presentation for
the cohomology of Grassmannians [DZ$_{1}$], we extend the homology rigidity
results known for the classical Grassmanians to the exceptional cases.
\end{abstract}



\section{Introduction}

Let $G$ be a compact connected Lie group with Lie algebra $L(G)$,
exponential map $\exp :L(G)\rightarrow G$, a maximal torus $T$, and a set $%
\Omega =\{\omega _{{1}},\cdots ,\omega _{{n}}\}\subset L(T)$ of \textsl{%
fundamental dominant weights} of $G$ (cf.\textbf{\ }[Hu]). For a weight $%
\omega \in \Omega $ let $H$ be the centralizer of the $1$--parameter
subgroup $\{\exp (t\omega )\in G\mid t\in \mathbb{R}\}$. The homogeneous $%
G/H $ is a \textsl{flag variety}, called \textsl{the Grassmannian} \textsl{%
of }$G $\textsl{\ corresponding to }$\omega $.

If $G$ is simple with rank $n$, we assume that the set $\Omega =\{\omega _{{1%
}},\cdots ,\omega _{{n}}\}$ is so ordered as the root--vertices in the
Dynkin diagram of $G$ pictured in [Hu, p.58]. With this convention we
specify, in the table below, eight Grassmanianns $X=G/H$ associated to the
exceptional Lie groups.

\begin{center}
{\normalsize 
\begin{tabular}{|c||c|c|c|c|c|c|c|c|}
\hline\hline
$G$ & $F_{4}$ & $F_{4}$ & $E_{6}$ & $E_{6}$ & $E_{7}$ & $E_{7}$ & $E_{7}$ & $%
E_{8}$ \\ \hline
$\omega $ & $\omega _{1}$ & $\omega _{4}$ & $\omega _{2}$ & $\omega _{6}$ & $%
\omega _{1}$ & $\omega _{7}$ & $\omega _{2}$ & $\omega _{8}$ \\ 
$H$ & $C_{3}\cdot S^{1}$ & $B_{3}\cdot S^{1}$ & $A_{6}\cdot S^{1}$ & $%
D_{5}\cdot S^{1}$ & $D_{6}\cdot S^{1}$ & $E_{6}\cdot S^{1}$ & $A_{7}\cdot
S^{1}$ & $E_{7}\cdot S^{1}$ \\ \hline
\end{tabular}
} \\[0pt]
\end{center}

For a topological space $X$ let $[X,X]$ be the set of homotopy classes of
self-maps of $X$, and let $End(H^{\ast }(X))$ be the set of endomorphisms of
the integral cohomology ring $H^{\ast }(X)$. It is known since 1970's that
if $X$ is a flag variety, the correspondence

\begin{center}
$[X,X]\rightarrow End(H^{\ast }(X))$
\end{center}

\noindent assigning a map $f$ with its induced endomorphism $f^{\ast }$ is
\textquotedblleft nearly faithful\textquotedblright\ in the sense that it
has finite kernel and finite cokernel with respect to the obvious monoid
structure on both sets (cf. [GH]). Therefore, the problem of finding $%
End(H^{\ast }(X))$ for a flag manifold $X$ is a crucial step toward
understanding all self-maps of $X$.

In this paper, we determine $End(H^{\ast }(X))$ for the eight Grassmannians
tabulated above. Before stating our result, it is worthwhile to mention two
earlier results on this topic.

\bigskip

It is known that, if $X$ is a Grassmannian, then

\begin{enumerate}
\item[(1.1)] $H^{r}(X)$ is torsion free and vanishes for odd $r$;

\item[(1.2)] $H^{2}(X)=\mathbb{Z}$ is generated by the Kaehlerian class $%
\alpha _{X}$ of $X$.
\end{enumerate}

\noindent We note that, in view of (1.2), any endomorphism of $H^{\ast }(X)$
must map $\alpha _{X}$ to a multiple of itself.

Let $A=\underset{r\geq 0}{\oplus }A^{2r}$ be a graded ring which vanishes in
odd degrees (as being indicated), and let $p$ be an integer. The element $%
\psi ^{p}\in End(A)$ defined by

\begin{quote}
$\psi ^{p}(x)=p^{r}x$, $x\in A^{2r}$
\end{quote}

\noindent is called the \textbf{Adams operator} on $A$ with order $p$.

\bigskip

\noindent \textbf{Theorem 1 (Hoffmann [H, 1984])}. \textsl{Let }$X$\textsl{\
be the Grassmannian }$G_{n,k}$\textsl{\ of }$k$\textsl{-dimensional linear
subspaces in }$\mathbb{C}^{n}$\textsl{, and let }$h\in End(H^{\ast }(X))$%
\textsl{\ with }$h(\alpha _{X})=p\alpha _{X},$\textsl{\ }$p\neq 0$\textsl{.
Then}

\textsl{i) if }$n\neq 2k$\textsl{, }$h=\psi ^{p}$\textsl{;}

\textsl{ii) if }$n=2k$\textsl{, there is the additional possibility }$h=\psi
^{p}\circ \tau $\textsl{, where }$\tau $\textsl{\ is the involution on }$%
G_{2k,k}$\textsl{\ that sends a }$k$\textsl{-dimensional linear subspaces }$%
L\subset C^{n}$\textsl{\ to its othogonal complement }$L^{\perp }$\textsl{.}$%
\square $

\bigskip

\noindent T\textbf{heorem 2 (Duan, [D, 2003]).} \textsl{Let }$X$\textsl{\ be
the Grassmannian }$\mathbb{C}S_{n}$\textsl{\ of complex structures on }$%
R^{2n}$\textsl{\ and let }$h\in End(H^{\ast }(X))$\textsl{\ with }$\
h(\alpha _{X})=p\alpha _{X}, $\textsl{\ }$p\neq 0$\textsl{. Then }$h=\psi
^{p}$\textsl{.}$\square $

\bigskip

Hoffmann conjectured that Theorem 1 holds for the case $p=0$. In comparison
it was shown in [D] that Theorem 2 is not true for $p=0$.

Turning to the exceptional Grassmannians under consideration, the case $%
X=E_{6}/A_{6}\cdot S^{1}$ will be of particular interest. Fix a maximal
torus $T$ in $E_{6}$. Then

\begin{quote}
a) the Dynkind diagram of $E_{6}$ has the symmetry $\varphi $ given by
reflection in the edge joinning the $2^{nd}$ and $4^{th}$ vertices [Hu,
p.58], hence determines an automorphism $\varphi ^{\ast }$ on $H^{\ast
}(E_{6}/T)$ by [P];

b) the induced ring map $p^{\ast }:$ $H^{\ast }(E_{6}/A_{6}\cdot
S^{1})\rightarrow H^{\ast }(E_{6}/T)$ of the fibration $p:E_{6}/T\rightarrow
E_{6}/A_{6}\cdot S^{1}$ (corresponding to the inclusion $T\subset A_{6}\cdot
S^{1}$) is injective, hence identify $H^{\ast }(E_{6}/A_{6}\cdot S^{1})$ as
a subring of $H^{\ast }(E_{6}/T)$.
\end{quote}

\noindent The main result of this paper is

\bigskip

\noindent \textbf{Theorem 3}. \textsl{Let }$X$\textsl{\ be one of the
Grassmanians specified in the Table, and let }$h\in End(H^{\ast }(X))$ 
\textsl{with} $h(\alpha _{X})=p\alpha _{X}$\textsl{. Then }

\textsl{i)} \textsl{if }$X\neq E_{6}/A_{6}\cdot S^{1},$\textsl{\ }$h=\psi
^{p}$\textsl{;}

\textsl{ii)} \textsl{if }$X=E_{6}/A_{6}\cdot S^{1}$\textsl{, there is the
additional possibility }$h=\psi ^{p}\circ \tau $\textsl{, where }$\tau $%
\textsl{\ is the restriction of }$\varphi ^{\ast }\in Aut(H^{\ast
}(E_{6}/T)) $ \textsl{to the subring }$H^{\ast }(E_{6}/A_{6}\cdot S^{1})$%
\textsl{.}

\bigskip

By establishing Theorem 3, we wish to demonstrate how the theory of Gr\"{o}%
bner basis is used in extending the stricking rigidity property to
Grassmannianns associated to the exceptional Lie groups.

The authors would like to thank Xu-an Zhao for clearifying the additional
case ii) in Theorem 3, which first occurs to us from the calulation in Case
3 of \S 4.

\section{The ring $H^{\ast}(X;\mathbb{Q})$ for a Grassmannian $X$}

If $X$ is a Grassmannian, the inclusion $\mathbb{Z}\rightarrow\mathbb{Q}$ of
coefficients induces an injective ring map $H^{\ast}(X)\rightarrow H^{\ast
}(X;\mathbb{Q})$ by (2.1). Therefore, it suffices to establish the Theorem
for cohomology with rational coefficients (instead of integers).

In the classical approaches to $End(H^{\ast }(X;\mathbb{Q}))$ for $X=G_{n,k%
\text{ }}$, $\mathbb{C}S_{n}$ (i.e. Theorems 1 and 2), the authors based
their calculation on the existing presentations of $H^{\ast }(X;\mathbb{Q})$
in the form of a quotient of a free polynomial ring

\begin{enumerate}
\item[(2.1)] $H^{\ast }(X;\mathbb{Q})=\mathbb{Q}[y_{1},\cdots
,y_{n}]/\left\langle g_{1},\cdots ,g_{m}\right\rangle $,
\end{enumerate}

\noindent where $g_{1},\cdots ,g_{m}\in $ $\mathbb{Q}[y_{1},\cdots ,y_{n}]$,
and where $\left\langle g_{1},\cdots ,g_{m}\right\rangle $ is the ideal
generated by $g_{1},\cdots ,g_{m}$. In [DZ$_{1}$], a unified method
computing the rings $H^{\ast }(X)$ for all Grassmannians $X$ has been
developed. In particular, for those Grassmannians $X$ concerned by the
Theorem, presentations of $H^{\ast }(X;\mathbb{Q})$ in the form of (2.1) can
be deduced directly from [DZ$_{1}$, Theorems 1--7] and [DZ$_{2}$; (2.9)]. We
present these results in (2.2)--(2.9) below, where the subscripts for the
generators $y_{i}$'s and relations $g_{j}$'s are adopted to indicate their
degrees in the fashion

\begin{quote}
$\deg y_{i}=2i$, $\deg g_{j}=2j$;
\end{quote}

\noindent and where $y_{1}$ is in the place of the Kaehlerian class $\alpha
_{X}$ on $X$ (see (1.2)).

\begin{enumerate}
\item[(2.2)] $H^{\ast}(F_{4}/C_{3}\cdot S^{1};\mathbb{Q})=\mathbb{Q}%
[y_{1},y_{4}]/\left\langle g_{8},g_{12}\right\rangle $, where

$g_{8}=24y_{4}^{2}+y_{1}^{8}-12y_{1}^{4}y_{4};~~$

$g_{12}=y_{1}^{12}-24y_{1}^{8}y_{4}+144y_{1}^{4}y_{4}^{2}-64y_{4}^{3}.$

\item[(2.3)] $H^{\ast}(F_{4}/B_{3}\cdot S^{1};\mathbb{Q})=\mathbb{Q}%
[y_{1},y_{4}]/\left\langle g_{8},g_{12}\right\rangle $, where

$g_{8}=3y_{4}^{2}-y_{1}^{8};~~$

$g_{12}=26y_{4}^{3}-5y_{1}^{12}.$

\item[(2.4)] $H^{\ast}(E_{6}/A_{6}\cdot S^{1};\mathbb{Q})=\mathbb{Q}%
[y_{1},y_{3},y_{4}]/\left\langle g_{8},g_{9},g_{12}\right\rangle $, where

$%
g_{8}=6y_{4}^{2}-12y_{1}y_{3}y_{4}+9y_{1}^{2}y_{3}^{2}+3y_{1}^{4}y_{4}-6y_{1}^{5}y_{3}+y_{1}^{8}; 
$

$%
g_{9}=-2y_{3}^{3}+6y_{3}y_{1}^{2}y_{4}-3y_{1}^{3}y_{3}^{2}+4y_{3}y_{1}^{6}-3y_{1}^{5}y_{4}-y_{1}^{9}; 
$

$%
g_{12}=4y_{4}^{3}-y_{3}^{4}+6y_{3}^{2}y_{1}^{2}y_{4}-4y_{3}^{3}y_{1}^{3}-2y_{3}^{2}y_{1}^{6}-9y_{1}^{4}y_{4}^{2}+12y_{1}^{5}y_{4}y_{3} 
$

$\qquad -6y_{1}^{8}y_{4}+4y_{1}^{9}y_{3}-y_{1}^{12}.$

\item[(2.5)] $H^{\ast }(E_{6}/D_{5}\cdot S^{1};\mathbb{Q})=\mathbb{Q}%
[y_{1},y_{4}]/\left\langle g_{9},g_{12}\right\rangle $, where

$g_{9}=2y_{1}^{9}+3y_{1}y_{4}^{2}-6y_{1}^{5}y_{4};$

$g_{12}=y_{4}^{3}-6y_{1}^{4}y_{4}^{2}+y_{1}^{12}.$

\item[(2.6)] $H^{\ast }(E_{7}/E_{6}\cdot S^{1};\mathbb{Q})=\mathbb{Q}%
[y_{1},y_{5},y_{9}]/\left\langle g_{10},g_{14},g_{18}\right\rangle $, where

$g_{10}=y_{5}^{2}-2y_{1}y_{9};$

$g_{14}=2y_{5}y_{9}-9y_{1}^{4}y_{5}^{2}+6y_{1}^{9}y_{5}-y_{1}^{14};$

$g_{18}=y_{9}^{2}+10y_{1}^{3}y_{5}^{3}-9y_{1}^{8}y_{5}^{2}+2y_{1}^{13}y_{5}.$

\item[(2.7)] $H^{\ast }(E_{7}/D_{6}\cdot S^{1};\mathbb{Q})=\mathbb{Q}%
[y_{1},y_{4},y_{6}]/\left\langle g_{12},g_{14},g_{18}\right\rangle $, where

$%
g_{12}=3y_{6}^{2}-y_{4}^{3}-3y_{1}^{4}y_{4}^{2}-2y_{1}^{6}y_{6}+2y_{1}^{8}y_{4}; 
$

$%
g_{14}=3y_{4}^{2}y_{6}+3y_{1}^{2}y_{6}^{2}+6y_{1}^{2}y_{4}^{3}+6y_{1}^{4}y_{4}y_{6}-3y_{1}^{6}y_{4}^{2}-4y_{1}^{8}y_{6}-2y_{1}^{10}y_{4}+y_{1}^{14}; 
$

$%
g_{18}=45y_{4}^{4}y_{1}^{2}+120y_{4}^{2}y_{1}^{4}y_{6}+60y_{4}^{3}y_{1}^{6}-52y_{4}^{2}y_{1}^{10}-16y_{1}^{6}y_{6}^{2}+80y_{1}^{8}y_{6}y_{4} 
$

$\qquad
-96y_{1}^{12}y_{6}-48y_{1}^{14}y_{4}+28y_{1}^{18}+116y_{6}^{3}+180y_{1}^{2}y_{4}y_{6}^{2}. 
$

\item[(2.8)] $H^{\ast }(E_{8}/E_{7}\cdot S^{1};\mathbb{Q})=\mathbb{Q}%
[y_{1},y_{6},y_{10}]/\left\langle g_{20},g_{24},g_{30}\right\rangle ,$ where

$%
g_{20}=3y_{10}^{2}+10y_{1}^{2}y_{6}^{3}+18y_{1}^{4}y_{6}y_{10}-12y_{1}^{10}y_{10}-18y_{1}^{8}y_{6}^{2}+9y_{1}^{14}y_{6}-y_{1}^{20}; 
$

$%
g_{24}=5y_{6}^{4}+30y_{1}^{2}y_{6}^{2}y_{10}+15y_{1}^{4}y_{10}^{2}-15y_{1}^{14}y_{10}-15y_{1}^{12}y_{6}^{2}+10y_{1}^{18}y_{6}-y_{1}^{24}; 
$

$%
g_{30}=12y_{1}^{4}y_{6}y_{10}^{2}+24y_{1}^{14}y_{6}y_{10}+56y_{1}^{24}y_{6}-36y_{1}^{10}y_{10}^{2}-32y_{10}^{3}+4y_{6}^{5}-9y_{1}^{30} 
$

$%
\qquad-48y_{1}^{20}y_{10}+60y_{1}^{6}y_{6}^{4}-64y_{1}^{18}y_{6}^{2}+96y_{1}^{8}y_{10}y_{6}^{2}-44y_{1}^{12}y_{6}^{3}. 
$

\item[(2.9)] $H^{\ast}(E_{7}/A_{7}\cdot S^{1},\mathbb{Q})=\mathbb{Q}%
[y_{1},y_{3},y_{4},y_{5},y_{7}]/\left\langle
g_{8},g_{10},g_{12},g_{14},g_{18}\right\rangle $, where

$%
g_{8}=6y_{4}^{2}-4y_{3}y_{5}+4y_{1}y_{7}-12y_{1}y_{3}y_{4}+9y_{1}^{2}y_{3}^{2}+2y_{1}^{3}y_{5}+3y_{1}^{4}y_{4}-6y_{1}^{5}y_{3}+y_{1}^{8}; 
$

$g_{10}=y_{5}^{2}-2y_{3}y_{7}+y_{1}^{3}y_{7};$

$%
g_{12}=-4y_{4}^{3}-2y_{3}^{2}y_{1}^{6}+9y_{1}^{4}y_{4}^{2}+6y_{1}^{8}y_{4}-4y_{1}^{9}y_{3}-24y_{1}y_{4}y_{7}-12y_{1}^{5}y_{4}y_{3}+8y_{5}y_{7} 
$

$%
\qquad+4y_{3}y_{1}y_{4}^{2}-18y_{3}^{2}y_{1}^{2}y_{4}+12y_{3}y_{1}^{2}y_{7}-3y_{3}^{4}+y_{1}^{12}+16y_{3}^{3}y_{1}^{3}; 
$

$%
g_{14}=y_{7}^{2}+3y_{4}y_{5}^{2}+y_{5}y_{3}^{3}+3y_{5}y_{3}y_{1}^{2}y_{4}-3y_{5}y_{1}^{3}y_{3}^{2}+y_{5}y_{3}y_{1}^{6}-y_{5}y_{1}y_{4}^{2}-3\,y_{5}y_{1}^{2}y_{7}; 
$

$%
g_{18}=-8y_{4}^{3}y_{5}y_{1}-4y_{5}y_{7}y_{3}^{2}-2y_{3}y_{1}^{7}y_{4}^{2}+y_{3}^{6}-6y_{3}y_{1}^{8}y_{7}+18y_{1}^{5}y_{3}^{2}y_{7}+15y_{3}^{2}y_{1}^{4}y_{4}^{2} 
$

$%
\qquad-6y_{3}y_{1}^{3}y_{4}^{3}+8y_{4}y_{5}y_{3}y_{1}^{6}-12y_{1}y_{5}^{2}y_{7}+y_{1}^{2}y_{4}^{4}+12y_{4}^{2}y_{5}^{2}+y_{3}^{2}y_{1}^{12}-6y_{1}^{9}y_{3}^{3} 
$

$%
\qquad+11y_{3}^{4}y_{1}^{6}+8y_{4}y_{5}y_{3}^{3}+24y_{4}^{2}y_{5}y_{3}y_{1}^{2}+8y_{3}y_{5}^{3}-18y_{3}^{3}y_{1}^{5}y_{4}-2y_{3}^{3}y_{1}y_{4}^{2} 
$

$%
\qquad+4y_{5}y_{7}y_{1}^{6}+6y_{3}^{4}y_{1}^{2}y_{4}-6y_{3}^{5}y_{1}^{3}+6y_{1}^{3}y_{4}^{2}y_{7}-18y_{3}y_{1}^{4}y_{4}y_{7}-6y_{3}^{3}y_{1}^{2}y_{7} 
$

$%
\qquad+6y_{3}^{2}y_{1}^{8}y_{4}-4y_{4}y_{7}^{2}+9y_{1}^{4}y_{7}^{2}-24y_{4}y_{5}y_{1}^{3}y_{3}^{2}-8y_{5}y_{7}y_{1}^{3}y_{3}-12y_{4}y_{5}y_{1}^{2}y_{7}. 
$
\end{enumerate}

\noindent\textbf{Remark.} In (2.2)--(2.9), every generator $y_{i}$ is in
fact a Schubert class of degree $i$ on the corresponding $X$ [DZ$_{1}$].

\section{\protect\smallskip The algorithm\ }

As one can see from (2.2)--(2.9) that the rings $H^{\ast }(X;\mathbb{Q})$
vary considerably with respect to the Grassmannians $X$. However, with their
common feature (2.1) in mind we may create some notations useful in a
unified approach to $End(H^{\ast }(X;\mathbb{Q}))$.

For a polynomial ring $\mathbb{Q}[y_{{1}},\cdots,y_{{n}}]$ graded by $\deg
y_{{i}}>0$, its subspace of the homogeneous elements of degree $r$ has a
canonical basis

\begin{enumerate}
\item[(3.1)] $B^{r}=\{y^{\alpha}=y_{{1}}^{b_{{1}}}\cdots,y_{{n}}^{b{n}%
}\mid\alpha=(b_{{1}},\cdots,b_{{n}})\in\mathbb{N}^{n},\deg y^{\alpha}=r\}$,
\end{enumerate}

\noindent called the \textsl{monomial basis }of $\mathbb{Q}[y_{{1}},\cdots
,y_{{n}}]$ in degree $r$, where $\mathbb{N}^{n}$ is the set of all $n$%
--tuples $\alpha =(b_{{1}},\cdots ,b_{{n}})$ of non--negative integers,
considered as an ordered set with respect to the lexicographical order on $%
\mathbb{N}^{n}$.

For a subset $\{g_{1},\cdots,g_{m}\}\subset$ $\mathbb{Q}[y_{1},\cdots,y_{n}]$
let $\mathcal{G}$ be the Gr\"{o}bner basis of the ideal $\left\langle
g_{1},\cdots,g_{m}\right\rangle $ [E]. The package "\textsl{Gbasis}" in
MAPLE [E] has the function

\begin{quote}
\textsl{In:} A subset $\{g_{1},\cdots,g_{m}\}\subset$ $\mathbb{Q}%
[y_{1},\cdots,y_{n}]$;

\textsl{Out:} the Gr\"{o}bner basis $\mathcal{G}$ of $\left\langle
g_{1},\cdots,g_{m}\right\rangle $,
\end{quote}

\noindent whose importance is shown in the next result.

\bigskip

\noindent\textbf{Lemma 1.} Each polynomial $f\in\mathbb{Q}[y_{{1}},\cdots,y_{%
{n}}] $ determines a unique element $h\in\mathbb{Q}[y_{{1}},\cdots,y_{{n}}]$%
, called \textsl{the residue of }$f$\textsl{\ module} $\left\langle
g_{1},\cdots,g_{m}\right\rangle $, that satisfies

\begin{quote}
i) $f\equiv h$ $\func{mod}$ $\mathcal{G}$;

ii) $f\in \left\langle g_{1},\cdots ,g_{m}\right\rangle $ if and only if $%
h=0 $.
\end{quote}

\noindent Moreover, the package "\textsl{Normal}" in MAPLE [E] has the
function to implement $h$ from $f$.$\square $

\bigskip

We may now clarify the algorithm by which Theorem 3 is established. An
endomorphism $f$ of the ring (2.1) can be regarded as an endomorphism of the
free polynomial ring $\mathbb{Q}[y_{1},\cdots ,y_{n}]$ which preserves the
ideal $\left\langle g_{1},\cdots ,g_{m}\right\rangle $. Based on this
observation we can assume that $f$ is given by

\begin{enumerate}
\item[(3.2)] $f(y_{i})=\sum\limits_{y^{\alpha}\in
B^{s_{i}}}c_{\alpha}^{i}\cdot y^{\alpha},~s_{i}=\deg y_{i},~c_{\alpha}^{i}\in%
\mathbb{Q},~i=1,\cdots,n$,
\end{enumerate}

\noindent that are subject to the restrictions

\begin{enumerate}
\item[(3.3)] $f(g_{j}(y_{1},\cdots,y_{n}))=g_{j}(f(y_{1}),\cdots,f(y_{n}))%
\in\left\langle g_{1},\cdots,g_{m}\right\rangle $, $~j=1,\cdots,m$,
\end{enumerate}

\noindent where the equality in (3.3) comes from the fact that $f$ is a ring
map. Assume that the residue of $g_{j}(f(y_{1}),\cdots,f(y_{n}))$ module $%
\left\langle g_{1},\cdots,g_{m}\right\rangle $ is

\begin{enumerate}
\item[(3.4)] $h_{j}=\sum\limits_{y^{\alpha}\in B^{t_{j}}}b_{\alpha}^{j}\cdot
y^{\alpha}$, $t_{j}=\deg g_{j}$,
\end{enumerate}

\noindent where $b_{\alpha}^{j}$ are certain polynomials in the $c_{\alpha
}^{i}$'s in (3.2). Then (3.3) is equivalent to the system

\begin{enumerate}
\item[(3.5)] $\mathcal{S}:\{b_{\alpha}^{j}=0\mid1\leq j\leq m,$ $%
y^{\alpha}\in B^{t_{j}}\}$
\end{enumerate}

\noindent by the Lemma. In other word, (3.5) constitutes all the constraints
that the $c_{\alpha }^{i}$'s in (3.2) must satisfy. Finally, we remark that
the package "\textsl{Solve}" or "\textsl{Gsolve}" in MAPLE [E] has the
function to produce all the solutions to the system (3.5).

\section{Proof of the Theorem}

\noindent \textbf{Case 1.} $X=F_{4}/C_{3}\cdot S^{1}$ (see (2.2)). In
accordance with the order $y_{1}>y_{4}$, the Gr\"{o}bner basis of the ideal $%
\left\langle g_{8},g_{12}\right\rangle $ is:

\begin{quote}
$\mathcal{G}=\left\{
24y_{4}^{2}+y_{1}^{8}-12y_{1}^{4}y_{4},~3y_{1}^{4}y_{4}^{2}-28y_{4}^{3},y_{4}^{4}\right\} 
$.
\end{quote}

\noindent Assume that $f\in End(H^{\ast}(X;\mathbb{Q}))$ is given by

\begin{quote}
$f(y_{1})=ky_{1},~~f(y_{4})=ay_{1}^{4}+by_{4},~~~k,a,b\in\mathbb{Q}$.
\end{quote}

\noindent Then the residues of $f(g_{i}),$ $i=8,12$,\ module $\left\langle
g_{8},g_{12}\right\rangle $ are (Lemma 1)

\begin{quote}
$h_{8}=12(24a^{2}+k^{8}-12k^{4}a+4ba-k^{4}b)y_{1}^{4}y_{4}$

$\qquad+12(-48a^{2}+24k^{4}a-2k^{8}+2b^{2})y_{4}^{2}$

$h_{12}=-(-1344k^{4}{b}^{2}+1792b^{2}a-832k^{12}+53248a^{3}-25344k^{4}ba$

$%
\qquad+64b^{3}-119808k^{4}a^{2}+19968k^{8}a+2112k^{8}b+16896ba^{2})y_{4}^{3} 
$.
\end{quote}

\noindent Applying \textsl{Solve} to the corresponding system (3.5) yields
that

\begin{quote}
$a=0,b=k^{4},k=k$.
\end{quote}

\noindent \textbf{Case 2.} $X=F_{4}/B_{3}\cdot S^{1}$ (see (2.3)). In
accordance with the order $y_{1}>y_{4}$, the Gr\"{o}bner basis of the ideal $%
\left\langle g_{8},g_{12}\right\rangle $ is:

\begin{quote}
$\mathcal{G}=\left\{
y_{1}^{8}-3y_{4}^{2},~15y_{4}^{2}y_{1}^{4}-26y_{4}^{3},~y_{4}^{4}\right\} $.
\end{quote}

\noindent Assume that $f\in End(H^{\ast}(X;\mathbb{Q}))$ is given by

\begin{quote}
$f(y_{1})=ky_{1},~~f(y_{4})=ay_{1}^{4}+by_{4},~~~k,a,b\in\mathbb{Q}$.
\end{quote}

\noindent Then the residues of $f(g_{i}),$ $i=8,12$,\ module $\left\langle
g_{8},g_{12}\right\rangle $ are (Lemma 1)

\begin{quote}
$h_{8}=3(3a^{2}-k^{8}+b^{2})y_{4}^{2}+6aby_{1}^{4}y_{4};$

$h_{12}=\frac{1}{5}%
(676a^{3}-130k^{12}+676ab^{2}+1170a^{2}b+130b^{3})y_{4}^{3}$
\end{quote}

\noindent Applying \textsl{Solve} to the corresponding system (3.5) yields
that

\begin{quote}
$a=0,b=k^{4},k=k$.
\end{quote}

\noindent \textbf{Case 3.} $X=E_{6}/A_{6}\cdot S^{1}$ (see (2.4)). In
accordance with the order $y_{1}>y_{3}>y_{4}$, the Gr\"{o}bner basis of the
ideal $\left\langle g_{8},g_{9},g_{12}\right\rangle $ is

\begin{quote}
$\mathcal{G}%
=\{6y_{4}^{2}-12y_{1}y_{3}y_{4}+9y_{1}^{2}y_{3}^{2}+3y_{1}^{4}y_{4}-6y_{1}^{5}y_{3}; 
$

$\qquad
y_{3}y_{1}^{6}-3y_{1}^{3}y_{3}^{2}+3y_{3}y_{1}^{2}y_{4}+y_{3}^{3}-3y_{4}^{2}y_{1};~ 
$

$%
\qquad3y_{3}^{2}y_{1}^{5}-8y_{3}^{3}y_{1}^{2}-3y_{4}^{2}y_{1}^{3}+12y_{4}y_{3}^{2}y_{1}-6y_{4}^{2}y_{3}; 
$

$%
\qquad2y_{3}^{4}-30y_{3}^{2}y_{1}^{2}y_{4}+33y_{4}^{2}y_{3}y_{1}-y_{3}^{3}y_{1}^{3}-22y_{4}^{3}+12y_{1}^{5}y_{4}y_{3};\cdots\} 
$,
\end{quote}

\noindent where we have omitted the elements in $\mathcal{G}$ with degree $%
>\max\{\deg g_{i}\}$ since they play no further role in the process. Assume
that $f\in End(H^{\ast}(X;\mathbb{Q}))$ is given by

\begin{quote}
$%
f(y_{1})=ky_{1},~~f(y_{3})=a_{1}y_{1}^{3}+a_{2}y_{3},~~f(y_{4})=b_{1}y_{1}^{4}+b_{2}y_{1}y_{3}+b_{3}y_{4}, 
$
\end{quote}

\noindent where $k,a_{1},a_{2},b_{1},b_{2},b_{3}\in \mathbb{Q}$. Then the
residues of $f(g_{i}),$ $i=8,9$,\ module $\left\langle
g_{8},g_{9},g_{12}\right\rangle $ are (Lemma 1)

\begin{quote}
$%
h_{8}=3(4b_{2}b_{3}+24b_{1}^{2}+36k^{2}a_{1}^{2}-48ka_{1}b_{1}+12k^{4}b_{1}-24k^{5}a_{1}+4k^{8} 
$

$%
\qquad-4ka_{2}b_{3})y_{1}y_{3}y_{4}+3(12b_{1}^{2}+2k^{8}-24ka_{1}b_{1}+18k^{2}a_{1}^{2}+6k^{4}b_{1} 
$

$\qquad-12k^{5}a_{1}+4b_{1}b_{2}+{k}%
^{4}b_{2}-2k^{5}a_{2}-4ka_{1}b_{2}-4ka_{2}b_{1}+6k^{2}a_{1}a_{2})y_{1}^{5}y_{3} 
$

$%
\qquad+3(3k^{2}a_{2}^{2}-18b_{1}^{2}-4ka_{2}b_{2}-9k^{4}b_{1}-27k^{2}a_{1}^{2}+2b_{2}^{2}+36ka_{1}b_{1}-3k^{8} 
$

$%
\qquad+18k^{5}a_{1})y_{1}^{2}y_{3}^{2}+3(12k^{5}a_{1}-2k^{8}-6k^{4}b_{1}+2b_{3}^{2}-12b_{1}^{2}+24ka_{1}b_{1} 
$

$%
\qquad-18k^{2}a_{1}^{2})y_{4}^{2}+3(k^{4}b_{3}-k^{8}+6k^{5}a_{1}-3k^{4}b_{1}-4ka_{1}b_{3}-6b_{1}^{2}+4b_{1}b_{3} 
$

$\qquad+12ka_{1}b_{1}-9k^{2}a_{1}^{2})y_{1}^{4}y_{4};$

$%
h_{9}=-18a_{1}^{2}a_{2}+12k^{6}a_{2}+72k^{2}a_{1}b_{1}-24a_{1}^{3}-12k^{9}-36k^{3}a_{1}^{2}+48k^{6}a_{1} 
$

$%
\qquad-36k^{5}b_{1}+18k^{2}a_{1}b_{2}+18k^{2}a_{2}b_{1}-18k^{3}a_{1}a_{2}-9k^{5}b_{2})y_{4}^{2}y_{1} 
$

$%
\qquad+(36k^{6}a_{1}-27k^{5}b_{1}-9k^{9}-27k^{3}a_{1}^{2}-18a_{1}^{2}a_{2}+12k^{6}a_{2}-9k^{5}b_{2} 
$

$%
\qquad-6a_{1}a_{2}^{2}+18k^{2}a_{1}b_{2}+18k^{2}a_{2}b_{1}-18k^{3}a_{1}a_{2}+54k^{2}a_{1}b_{1}-18a_{1}^{3} 
$

$%
\qquad+6k^{2}a_{2}b_{2}-3k^{3}a_{2}^{2})y_{1}^{3}y_{3}^{2}+(9k^{5}b_{1}+6a_{1}^{3}+3k^{9}+9k^{3}a_{1}^{2}-12k^{6}a_{1} 
$

$%
\qquad-18k^{2}a_{1}b_{1}+6k^{2}a_{1}b_{3}-3k^{5}b_{3})y_{1}^{5}y_{4}+(-24k^{6}a_{1}+18k^{5}b_{1}+6a_{1}^{2}a_{2} 
$

$%
\qquad-4k^{6}a_{2}+3k^{5}b_{2}+6k^{9}-36k^{2}a_{1}b_{1}-6k^{2}a_{1}b_{2}-6k^{2}a_{2}b_{1}+6k^{3}a_{1}a_{2} 
$

$%
\qquad+12a_{1}^{3}+18k^{3}a_{1}^{2}-2a_{2}^{3})y_{3}^{3}+(6k^{9}+12a_{1}^{3}+18a_{1}^{2}a_{2}-18k^{2}a_{1}b_{2} 
$

$%
\qquad-18k^{2}a_{2}b_{1}+18k^{3}a_{1}a_{2}-12k^{6}a_{2}+9k^{5}b_{2}+18k^{3}a_{1}^{2}-24k^{6}a_{1} 
$

$\qquad+18k^{5}b_{1}-36k^{2}a_{1}b_{1}+6k^{2}a_{2}b_{3})y_{3}y_{1}^{2}y_{4}$.
\end{quote}

\noindent With $k,a_{1},a_{2},b_{1},b_{2},b_{3}\in\mathbb{Q}$ the
corresponding system (3.5) has two solutions

\begin{quote}
$\{a_{1}=k^{3},a_{2}=-k^{3},b_{1}=k^{4},b_{2}=-2k^{4},b_{3}=k^{4}\},$

$\{a_{1}=0,a_{2}=k^{3},b_{1}=0,b_{2}=0,b_{3}=k^{4}\}$,
\end{quote}

\noindent and with respect to the first solution, the residue $h_{12}$ of $%
f(g_{12})$ is trivial. This establish the Theorem for the current situation.

\bigskip

\noindent \textbf{Case 4.} $X=E_{6}/D_{5}\cdot S^{1}$ (see (2.5)). In
accordance with the order $y_{1}>y_{4}$, the Gr\"{o}bner basis of the ideal $%
\left\langle g_{9},g_{12}\right\rangle $ is:

\begin{quote}
$\mathcal{G}%
=\{2y_{1}^{9}+3y_{1}y_{4}^{2}-6y_{1}^{5}y_{4},~6y_{4}y_{1}^{8}-15y_{1}^{4}y_{4}^{2}+2y_{4}^{3}, 
$ $3y_{4}^{2}y_{1}^{5}-7y_{4}^{3}y_{1},$

$\qquad y_{4}^{3}y_{1}^{4}-2y_{4}^{4},~y_{4}^{4}y_{1},~y_{4}^{5}\}$.
\end{quote}

\noindent Assume that $f\in End(H^{\ast}(X;\mathbb{Q}))$ is given by

\begin{quote}
$f(y_{1})=ky_{1},~~f(y_{4})=ay_{1}^{4}+by_{4},~~~k,a,b\in\mathbb{Q}$.
\end{quote}

\noindent Then the residues of $f(g_{i}),$ $i=9,12$,\ module $\left\langle
g_{9},g_{12}\right\rangle $ are (Lemma 1)

\begin{quote}
$h_{9}=\frac{3}{2}(-2k^{8}+6k^{4}a-3a^{2}+2b^{2})ky_{1}y_{4}^{2}+\frac{3}{2}%
(4k^{8}-12k^{4}a+6a^{2}$

$\qquad-4k^{4}b+4ab)ky_{1}^{5}y_{4};$

$h_{12}=\frac{1}{2}%
(-2a^{3}-2k^{12}+12k^{4}a^{2}-2a^{2}b+2b^{3}+8k^{4}ab)y_{4}^{3}+\frac{1}{2}%
(12a^{3}$

$%
\qquad+12k^{12}-60k^{4}ab+6ab^{2}-12k^{4}b^{2}-72k^{4}a^{2}+15a^{2}b)y_{1}^{4}y_{4}^{2} 
$.
\end{quote}

\noindent Applying \textsl{Solve} to the corresponding system (3.5) yields
that

\begin{quote}
$a=0,b=k^{4},k=k$.
\end{quote}

\noindent \textbf{Case 5.} $X=E_{7}/E_{6}\cdot S^{1}$ (see (2.6)). In
accordance with the order $y_{1}>y_{4}>y_{9}$, the Gr\"{o}bner basis of the
ideal $\left\langle g_{10},g_{14},g_{18}\right\rangle $ is:

\begin{quote}
$\mathcal{G}=\{-2y_{5}y_{9}+9y_{1}^{4}y_{5}^{2}-6y_{1}^{9}y_{5}+y_{1}^{14}$, 
$y_{9}^{2}+10y_{1}^{3}y_{5}^{3}-9y_{1}^{8}y_{5}^{2}+2y_{1}^{13}y_{5},$

$\qquad-y_{5}^{2}+2y_{1}y_{9},\cdots\}$,
\end{quote}

\noindent where the elements in $\mathcal{G}$ with degree $>\max\{\deg
g_{i}\}$ has been omitted because they play no role in the latter course.
Assume that $f\in End(H^{\ast}(X;\mathbb{Q}))$ is given by

\begin{center}
$\qquad
f(y_{1})=ky_{1},~f(y_{5})=a_{1}y_{1}^{5}+a_{2}y_{5},~f(y_{9})=b_{1}y_{1}^{9}+b_{2}y_{1}^{4}y_{5}+b_{3}y_{9}, 
$
\end{center}

\noindent where $k,a_{1},a_{2},b_{1},b_{2},b_{3}\in \mathbb{Q}$. Then the
residues of $f(g_{i}),$ $i=10,14,18$, module $\left\langle
g_{10},g_{14},g_{18}\right\rangle $ are (Lemma 1)

\begin{quote}
$%
h_{10}=(2a_{2}^{2}-2kb_{3})y_{1}y_{9}+(2a_{1}a_{2}-2kb_{2})y_{1}^{5}y_{5}+(a_{1}^{2}-2kb_{1})y_{1}^{10}; 
$

$%
h_{14}=2(81k^{4}a_{1}^{2}+2a_{2}b_{2}-9k^{4}a_{2}^{2}-18a_{1}b_{1}-54k^{9}a_{1}+9k^{14}+a_{1}b_{3})y_{9}y_{1}^{5} 
$

$%
\qquad+2(a_{1}b_{2}+a_{2}b_{1}+3k^{9}a_{2}+6a_{1}b_{1}-27k^{4}a_{1}^{2}+18k^{9}a_{1}-9k^{4}a_{1}a_{2} 
$

$%
\qquad-3k^{14})y_{1}^{9}y_{5}+2(a_{2}b_{3}-k^{14}+2a_{1}b_{1}-9k^{4}a_{1}^{2}+6k^{9}a_{1})y_{5}y_{9}; 
$

$%
h_{18}=(2b_{1}b_{3}+270k^{3}a_{1}^{2}a_{2}-162k^{8}a_{1}a_{2}+72k^{13}a_{1}+18b_{1}b_{2}+18k^{13}a_{2} 
$

$%
\qquad-18k^{8}a_{2}^{2}+360k^{3}a_{1}^{3}-324k^{8}a_{1}^{2}+60k^{3}a_{1}a_{2}^{2}+2b_{2}^{2}+36b_{1}^{2})y_{9}y_{1}^{9} 
$

$%
\qquad+(522k^{8}a_{1}^{2}-580k^{3}a_{1}^{3}+2b_{2}b_{3}-300k^{3}a_{1}^{2}a_{2}-58b_{1}^{2}-116k^{13}a_{1} 
$

$%
\qquad+180k^{8}a_{1}a_{2}-20b_{1}b_{2}-20k^{13}a_{2}+20k^{3}a_{2}^{3})y_{9}y_{5}y_{1}^{4}+(b_{3}^{2}-15k^{3}a_{1}^{2}a_{2} 
$

$%
\qquad+9k^{8}a_{1}a_{2}-3b_{1}^{2}-k^{13}a_{2}-30k^{3}a_{1}^{3}+27k^{8}a_{1}^{2}-6k^{13}a_{1}-b_{1}b_{2})y_{9}^{2} 
$.
\end{quote}

\noindent Applying \textsl{Solve} to the corresponding system (3.5) yields
that

\begin{quote}
$a_{1}=b_{1}=b_{2}=0,a_{2}=k^{5},b=k^{9},k=k$.
\end{quote}

\noindent \textbf{Case 6.} $X=E_{7}/D_{6}\cdot S^{1}$, $E_{8}/E_{7}\cdot
S^{1}$ or $E_{7}/A_{7}\cdot S^{1}$. Procedure the same as the above
establish Theorem 3 for these cases. Only the expressions of the
corresponding $h_{i}$'s are very lengthy, and the package "\textsl{Gsolve}"
in MAPLE appears more effective than the "\textsl{Solve}" when implementing
the corresponding system (3.5). We omit the details and refer to the thesis
[L] of the first author for thorough discussions on these cases.

\begin{center}
\textbf{References}
\end{center}

\begin{enumerate}
\item[{[D]}] H. Duan, Self-maps of the Grassmannian of complex structures,
Compositio Math., 132 (2002), 159-175.

\item[{[DZ$_{1}$]}] H. Duan and Xuezhi Zhao, The Chow rings of generalized
Grassmannians, arXiv: math.AG/0511332.

\item[{[DZ$_{2}$]}] H. Duan and Xuezhi Zhao, The integral cohomology of
complete flag manifolds, arXiv: math.AT/0801.2444.

\item[{[E]}] K. Ernic, \textsl{A Guide To Maple}, Springer Verlag, 1999.

\item[{[GH]}] H. Glover and W. Homer, Self-maps of flag manifolds, Trans. AMS
267(1981), 423-434.

\item[{[H$_{1}$]}] M. Hoffman, Endomorphisms of the cohomology of complex
Grassmannians, Trans, AMS 281 (1984), 745-740.

\item[{[H$_{2}$]}] M. Hoffman, On fixed point free maps of the complex flag
manifold, Indiana Math. J., 33(1984), 249-255.

\item[{[Hu]}] J. E. Humphreys, \textsl{Introduction to Lie algebras and
representation theory,} Graduated Texts in Math. 9, Springer-Verlag New
York, 1972.

\item[{[L]}] F. Li, Endomorphisms of the cohomology ring of a generalized
Grassmannian, Thesis, Jilin University.

\item[{[P]}] S. Papadima, Rigidity properties of compact Lie groups modulo
maximal tori, Math. Ann. 275(1987), 637--652.
\end{enumerate}

\end{document}